\documentclass{elsart}

\usepackage{amssymb}
\usepackage{amsmath}
\usepackage{lscape}

\numberwithin{equation}{section}
\numberwithin{thm}{section}
\numberwithin{table}{section}

\newcommand{\ft}{\Psi}

\newcommand{\s}{\varphi}

\newcommand{\sfrak}{\mathfrak{s}}
\newcommand{\g}{\mathfrak{g}}

\allowdisplaybreaks

\begin{document}

\begin{frontmatter}
\journal{J. Math. Anal. Appl.}

\title{Ramanujan-Slater Type Identities \\Related to the Moduli 18 and 24}
\date{February 21, 2008}

\author{James McLaughlin}
\address{Department of Mathematics,
West Chester University,
West Chester, PA;
telephone 610-738-0585; fax 610-738-0578}
\ead{jmclaughl@wcupa.edu}
\ead[url]{http://math.wcupa.edu/\~{}mclaughlin}

\author{Andrew V. Sills}
\address{Department of Mathematical Sciences,
Georgia Southern University,
Statesboro, GA;
telephone 912-681-5892;  fax 912-681-0654}
\ead{asills@GeorgiaSouthern.edu}
\ead[url]{http://math.georgiasouthern.edu/\~{}asills}

\begin{abstract} We present several new families of Rogers-Ramanujan type identities
related to the moduli 18 and 24.  A few of the identities were found by
either Ramanujan, Slater,
or Dyson, but most are believed to be new.
For one of these families, we discuss possible connections with Lie algebras.
We also present two families of related
false theta function identities.
\end{abstract}

\end{frontmatter}

\section{Introduction}
The Rogers-Ramanujan identities are
\begin{thm}[The Rogers-Ramanujan Identities]
 \begin{equation}\label{RRa1}
   \sum_{n=0}^\infty \frac{q^{n^2}}{(q;q)_n} =
   \frac{(q^2, q^3, q^5; q^5)_\infty}{(q;q)_\infty},
 \end{equation} and
\begin{equation}\label{RRa2}
   \sum_{n=0}^\infty \frac{q^{n(n+1)}}{(q;q)_n} =
   \frac{(q, q^4, q^5; q^5)_\infty}{(q;q)_\infty},
 \end{equation}
where \[ (a;q)_m = \prod_{j=0}^{m-1} (1-aq^j), \quad
      (a;q)_\infty = \prod_{j=0}^\infty (1-aq^j), \] and
      \[ (a_1, a_2, \dots, a_r; q)_s = (a_1;q)_s (a_2;q)_s \dots (a_r;q)_s. \]
\end{thm}
(Although the results in this paper may be considered purely from the
point of view of formal power series, they also yield identities of
analytic functions provided $|q|<1$.)

The Rogers-Ramanujan identities are due to L.~J.~Rogers~\cite{R94},
and were rediscovered by S. Ramanujan~\cite{M18} and
I. Schur~\cite{S17}.
Rogers and others discovered many series--product
identities similar in form to the Rogers-Ramanujan identities, and
such identities are called ``identities of the Rogers-Ramanujan type."
Two of the largest collections of Rogers-Ramanujan type identities are
contained in Slater's paper~\cite{S52} and Ramanujan's
Lost Notebook~\cite[Chapters 10--11]{AB05},~\cite[Chapters 1--5]{AB07}.

Rogers-Ramanujan type identities occur in closely related ``families."
Just as there are two Rogers-Ramanujan identities related to the modulus 5, there are
a family of three
Rogers-Selberg identities related to the modulus 7~\cite[p. 331, (6)]{R17}, a family of
three identities related to the
modulus 9 found by Bailey~\cite[p. 422, Eqs. (1.6)--(1.8)]{B47},
a family of four identities related to the modulus 27 found by
Dyson~\cite[p. 433, Eqs. (B1)--(B4)]{B47}, etc.

While both Ramanujan and Slater usually managed to find all members of
a given family, this was not always the case.  In this paper, we present
several complete
families of identities for which Ramanujan or Slater found only one member, as
well as two complete new families.

The following family of four identities related to the modulus 18
is believed to be new:
{\allowdisplaybreaks
\begin{gather}
\sum_{n=0}^\infty \frac{  q^{n(n+1)} (-1;q^3)_n}{ (-1;q)_n (q;q)_{2n} }
= \frac{(q,q^8,q^9;q^9)_\infty (q^7,q^{11};q^{18})_\infty}
{(q;q)_\infty} \label{m18-1}\\
\sum_{n=0}^\infty \frac{  q^{n^2} (-1;q^3)_n}{ (-1;q)_n (q;q)_{2n} }
= \frac{(q^2,q^7,q^9;q^9)_\infty (q^5,q^{13} ; q^{18})_\infty}
{(q;q)_\infty} \label{m18-2} \\
\sum_{n=0}^\infty \frac{  q^{n(n+1)} (-q^3;q^3)_n}{ (-q;q)_n (q;q)_{2n+1} }
= \frac{(q^3,q^6,q^9;q^9)_\infty (q^3,q^{15};q^{18})_\infty}{(q;q)_\infty}
\label{m18-3} \\
\sum_{n=0}^\infty \frac{  q^{n(n+2)} (-q^3;q^3)_n }
{ (q^2;q^2)_n (q^{n+2};q)_{n+1} }
= \frac{(q^4,q^5,q^9;q^9)_\infty (q,q^{17};q^{18})_\infty}{(q;q)_\infty}
\label{m18-4}
\end{gather}
}

\begin{rem}
We included Identity~\eqref{m18-3} in our joint paper with D. Bowman~\cite[Eq. (6.30)]{BMS07}, as
it also occurs as part of a different family of four identities.
\end{rem}

A closely related family of mod 18 identities is as follows.
\begin{gather}
1+\sum_{n=1}^\infty \frac{  q^{n^2} (q^3;q^3)_{n-1} (2+q^n)}
{ (q;q)_{n-1} (q;q)_{2n} }
= \frac{(-q,-q^8,q^9;q^9)_\infty (q^7,q^{11};q^{18})_\infty}{(q;q)_\infty}
\label{m18-m1}\\
1+\sum_{n=1}^\infty \frac{  q^{n^2} (q^3;q^3)_{n-1} (1+2q^n)}
{ (q;q)_{n-1} (q;q)_{2n} }
= \frac{(-q^2,-q^7,q^9;q^9)_\infty (q^5,q^{13};q^{18})_\infty}
{(q;q)_\infty} \label{m18-m2}\\
 \sum_{n=0}^\infty \frac{  q^{n(n+1)} (q^3;q^3)_n}{ (q;q)_n (q;q)_{2n+1} }
= \frac{(-q^3,-q^6,q^9;q^9)_\infty (q^3,q^{15};q^{18})_\infty}{(q;q)_\infty}
 \label{m18-m3}\\
\sum_{n=0}^\infty \frac{  q^{n(n+2)} (q^3;q^3)_n }
{ (q;q)_n^2 (q^{n+2};q)_{n+1} }
= \frac{(-q^4,-q^5,q^9;q^9)_\infty (q,q^{17};q^{18})_\infty}{(q;q)_\infty}
 \label{m18-m4}
\end{gather}
Identity~\eqref{m18-m3} is due to Dyson~\cite[p. 434, Eq. (B3)]{B47} and also appears in Slater~\cite[p. 161, Eq. (92)]{S52}.
In both~\cite{B47} and~\cite{S52}, the right hand side of~\eqref{m18-m3} appears
in a different form and thus is seen to be
a member of a different family of four identities related to the modulus 27.

Following Ramanujan (cf. \cite[p. 11, Eq (1.1.7)]{AB05}), let us use the notation
\begin{equation*}
\psi(q) = \frac{(q^2;q^2)_\infty}{(q;q^2)_\infty}.
\end{equation*}
Ramanujan recorded the identity
\begin{equation}
\sum_{n=0}^\infty \frac{  q^{n^2}  (-q^3;q^6)_n}{ (q^2;q^2)_{2n} }
= \frac{ (q^2,q^{10},q^{12};q^{12})_\infty (q^{8},q^{16};q^{24})_\infty }{\psi(-q)}  \label{m24t-2}
\end{equation}
in his lost notebook~\cite[Entry 5.3.8]{AB07}.
As we see below, it is actually only one of a family of five similar identities.
\begin{gather}
\sum_{n=0}^\infty \frac{  q^{n(n+2)} (-q;q^2)_n (-1;q^6)_n}{ (q^2;q^2)_{2n} (-1;q^2)_n}
= \frac{ (q,q^{11},q^{12};q^{12})_\infty (q^{10},q^{14};q^{24})_\infty }{\psi(-q)} \label{m24t-1}\\
\sum_{n=0}^\infty \frac{  q^{n^2}  (-q;q^2)_n (-1;q^6)_n}{ (q^2;q^2)_{2n} (-1;q^2)_n }
= \frac{ (q^3,q^{9},q^{12};q^{12})_\infty (q^{6},q^{18};q^{24})_\infty }{\psi(-q)} \label{m24t-3}\\
\sum_{n=0}^\infty \frac{  q^{n(n+2)}  (-q^3;q^6)_n}{ (q;q)_{2n+1} (-q;q)_{2n}  }
= \frac{ (q^4,q^{8},q^{12};q^{12})_\infty (q^{4},q^{20};q^{24})_\infty }{\psi(-q)} \label{m24t-4}\\
\sum_{n=0}^\infty \frac{  q^{n(n+2)}  (-q;q^2)_{n+1} (-q^6;q^6)_n }
{ (q^4;q^4)_{n} (q^{2n+4};q^2)_{n+1}  }
= \frac{ (q^5,q^{7},q^{12};q^{12})_\infty (q^{2},q^{22};q^{24})_\infty }{\psi(-q)} \label{m24t-5}
\end{gather}

Ramanujan also recorded the identity
\begin{equation} \label{m24t-m2}
\sum_{n=0}^\infty \frac{  q^{n^2}  (q^3;q^6)_n}{ (q;q^2)_{n}^2 (q^4;q^4)_n }
= \frac{ (-q^2,-q^{10},q^{12};q^{12})_\infty (q^{8},q^{16};q^{24})_\infty }{\psi(-q)}
\end{equation}
in the lost notebook~\cite[Entry 5.3.9]{AB07}.

Again, it is one of a family of five similar identities.  This time, however, two of the
remaining four identities were found by Slater.  Identity~\eqref{m24t-m4} is
a corrected presentation of~\cite[p. 164, Eq. (110)]{S52} and
identity~\eqref{m24t-m5} is
a corrected presentation of~\cite[p. 163, Eq. (108)]{S52}.

\begin{gather}
1+\sum_{n=1}^\infty \frac{  q^{n^2} (-q;q^2)_n (q^6;q^6)_{n-1} (2+q^{2n})}{ (q^2;q^2)_{2n} (q^2;q^2)_{n-1}}
= \frac{ (-q,-q^{11},q^{12};q^{12})_\infty (q^{10},q^{14};q^{24})_\infty }
{\psi(-q)} \label{m24t-m1}\\
1+\sum_{n=1}^\infty \frac{  q^{n^2}  (-q;q^2)_n (q^6;q^6)_{n-1} (1+2q^{2n})}
{ (q^2;q^2)_{2n} (q^2;q^2)_{n-1} }
= \frac{ (-q^3,-q^{9},q^{12};q^{12})_\infty (q^{6},q^{18};q^{24})_\infty }
{\psi(-q)} \label{m24t-m3}\\
 \sum_{n=0}^\infty \frac{  q^{n(n+2)}  (q^3;q^6)_n (-q;q^2)_{n+1} }{ (q^2;q^2)_{2n+1} (q;q^2)_n }
= \frac{ (-q^4,-q^{8},q^{12};q^{12})_\infty (q^{4},q^{20};q^{24})_\infty }
{\psi(-q)} \label{m24t-m4}\\
 \sum_{n=0}^\infty \frac{  q^{n(n+2)}  (-q;q^2)_{n+1} (q^6;q^6)_n }
{ (q^{2n+4};q^2)_{n+1} (q^2;q^2)_n^2 }
= \frac{ (-q^5,-q^{7},q^{12};q^{12})_\infty (q^{2},q^{22};q^{24})_\infty }
{\psi(-q)}\label{m24t-m5}
\end{gather}

We believe that the following family of five identities has not previously appeared
in the literature:
\begin{gather}
\sum_{n=0}^\infty \frac{  q^{n(n+1)} (-q^2;q^2)_n (-q^3;q^6)_n }
 { (q;q)_{2n} (-q;q)_{2n+1} (-q;q^2)_n}
= \frac{ (q,q^{11},q^{12};q^{12})_\infty (q^{10},q^{14};q^{24})_\infty }
{\varphi(-q^2)}\label{m24s-1}\\
\sum_{n=0}^\infty \frac{  q^{n(n+1)}  (-1;q^6)_n (-q^2;q^2)_n }{ (q^2;q^2)_{2n} (-1;q^2)_n }
= \frac{ (q^2,q^{10},q^{12};q^{12})_\infty (q^{8},q^{16};q^{24})_\infty }
{\varphi(-q^2)}
\label{m24s-2}\\
\sum_{n=0}^\infty \frac{  q^{n(n+1)}  (-q^2;q^2)_n (-q^3;q^6)_n}{ (q^2;q^2)_{2n+1} (-q;q^2)_n }
= \frac{ (q^3,q^{9},q^{12};q^{12})_\infty (q^{6},q^{18};q^{24})_\infty }
{\varphi(-q^2)}\label{m24s-3}\\
\sum_{n=0}^\infty \frac{  q^{n(n+1)}  (-q^6;q^6)_n}{ (q^2;q^2)_{2n+1}  }
= \frac{ (q^4,q^{8},q^{12};q^{12})_\infty (q^{4},q^{20};q^{24})_\infty }
{\varphi(-q^2)} \label{m24s-4}\\
\sum_{n=0}^\infty \frac{  q^{n(n+3)}  (-q^2;q^2)_{n} (-q^3;q^6)_n }
{ (q^2;q^2)_{2n+1} (-q;q^2)_n }
= \frac{ (q^5,q^{7},q^{12};q^{12})_\infty (q^{2},q^{22};q^{24})_\infty }
{\varphi(-q^2)} \label{m24s-5}
,\end{gather}
where
\begin{equation*}
  \varphi(q) := \frac{(-q;-q)_\infty}{(q;-q)_\infty}
 \end{equation*}
is another notation used by Ramanujan.

In the following counterpart to the preceding family, two of the five identities appear
in Slater's list.
\begin{gather}
\sum_{n=0}^\infty \frac{  q^{n(n+1)} (-q^2;q^2)_n (q^3;q^6)_{n} }{ (q;q)_{2n+1}
(-q;q)_{2n}
(q;q^2)_{n}}
= \frac{ (-q,-q^{11},q^{12};q^{12})_\infty (q^{10},q^{14};q^{24})_\infty }
{\varphi(-q^2)}\label{m24s-m1} \\
1+\sum_{n=1}^\infty \frac{  q^{n(n+1)}  (q^6;q^6)_{n-1} (-q^2;q^2)_n}{ (q^2;q^2)_{n-1} (q^2;q^2)_{2n} }
= \frac{ (-q^2,-q^{10},q^{12};q^{12})_\infty (q^{8},q^{16};q^{24})_\infty }
  {\varphi(-q^2)}  \label{m24s-m2}
\\
\sum_{n=0}^\infty \frac{  q^{n(n+1)}  (-q^2;q^2)_n (q^3;q^6)_{n}}{ (q^2;q^2)_{2n+1} (q;q^2)_{n} }
= \frac{ (-q^3,-q^{9},q^{12};q^{12})_\infty (q^{6},q^{18};q^{24})_\infty }{\varphi(-q^2)}\label{m24s-m3}\\
\sum_{n=0}^\infty \frac{  q^{n(n+1)}  (q^6;q^6)_n (-q^2;q^2)_n }
{ (q^2;q^2)_{2n+1} (q^2;q^2)_n }
= \frac{ (-q^4,-q^{8},q^{12};q^{12})_\infty (q^{4},q^{20};q^{24})_\infty }{\varphi(-q^2)}\label{m24s-m4}\\
 \sum_{n=0}^\infty \frac{  q^{n(n+3)}  (-q^2;q^2)_{n} (q^3;q^6)_n }
{ (q^2;q^2)_{2n+1} (q;q^2)_n }
= \frac{ (-q^5,-q^{7},q^{12};q^{12})_\infty (q^{2},q^{22};q^{24})_\infty} {\varphi(-q^2)}\label{m24s-m5}
\end{gather}
Identity~\eqref{m24s-m3} is due to Slater~\cite[p. 163, Eq. (107)]{S52}.
Identity~\eqref{m24s-m4} is originally due to
Dyson~\cite[p. 434, Eq. (D2)]{B47} and also appears in Slater~\cite[p. 160,
Eq. (77)]{S52}.

The following false theta series identities, which are closely
related to identities~\eqref{m24s-1}--\eqref{m24s-m5}, are believed
to be new, except for~\eqref{ft7} and ~\eqref{ft9}.  Identity~\eqref{ft7}
is due to Dyson~\cite[p. 434, Eq. (E1)]{B47}, while
Identity~\eqref{ft9} appears in Ramanujan's lost notebook~\cite[Entry 5.4.2]{AB07}
and was rediscovered by Dyson~\cite[p. 434, Eq. (E2)]{B47}.
\begin{multline}
\sum_{n=0}^\infty \frac{(-1)^n q^{n(n+1)} (-q^3;q^6)_n} { (q^{2};q^4)_n  (-q;q)_{2n+1} }\\
= \sum_{n=0}^\infty (-1)^n q^{18n^2 + 3n}(1+q^{30n+15}) - q \label{ft1} \sum_{n=0}^\infty
(-1)^n q^{18n^2 + 9n}(1+q^{18n+9})
\end{multline}
\begin{equation}
\sum_{n=0}^\infty \frac{ (-1)^n q^{n(n+3)} (-q^6;q^6)_n}{
(q^{2};q^4)_{n+1} (-q^2;q^2)_{n} (-q^2;q^2)_{n+1}}
=\sum_{n=0}^\infty (-1)^n q^{18n^2+12n} (1+q^{12n+6})  \label{ft2}
\end{equation}
\begin{multline}
\sum_{n=0}^\infty \frac{ (-1)^n q^{n(n+1)} (-q^3;q^6)_n}{
(q^{2};q^4)_{n+1}  (-q;q)_{2n}}
\\=\sum_{n=0}^\infty (-1)^n q^{18n^2+3n}(1+q^{30n+15})
  +q^3\sum_{n=0}^\infty (-1)^n q^{18n^2+15n}(1+q^{6n+3}) \label{ft3}
\end{multline}
\begin{multline}
\sum_{n=0}^\infty \frac{ (-1)^n q^{n(n+1)} (-q^6;q^6)_n}{ (q^{2};q^4)_{n+1} (-q^2;q^2)_{n}^2}\\
=\sum_{n=0}^\infty (-1)^n q^{18n^2+6n} (1+q^{24n+12})
 +2 q^4 \sum_{n=0}^\infty (-1)^n q^{18n^2+18n} \label{ft4}
\end{multline}
\begin{multline}
\sum_{n=0}^\infty \frac{ (-1)^n q^{n(n+3)} (-q^3;q^6)_n}{ (q^{2};q^4)_{n+1}  (-q;q)_{2n}}\\
=\sum_{n=0}^\infty (-1)^n q^{18n^2+9n} (1+q^{18n+9})
 + q^2 \sum_{n=0}^\infty (-1)^n q^{18n^2+15n} (1+q^{6n+3}) \label{ft5}
\end{multline}
\begin{multline}
\sum_{n=0}^\infty \frac{(-1)^n q^{n(n+1)} (q^3;q^6)_n }{
(q^{2};q^4)_n (-q^2;q^2)_n (q;q^2)_{n+1}  }
\\=\sum_{n=0}^\infty q^{18n^2 + 3n}(1-q^{30n+15})  + q \sum_{n=0}^\infty
q^{18n^2 + 9n}(1-q^{18n+9}) \label{ft6}
\end{multline}
\begin{equation}
\sum_{n=0}^\infty \frac{ (-1)^{n} q^{n(n+3)}(q^6;q^6)_{n}}{  (q;q)_{2n+1} (-q;q)_{2n+2} }
=\sum_{n=0}^\infty q^{18n^2+12n} (1-q^{12n+6})\label{ft7}
\end{equation}
\begin{multline}
\sum_{n=0}^\infty \frac{ (-1)^n q^{n(n+1)} (q^3;q^6)_n}{
(q^{2};q^4)_{n+1} (-q^2;q^2)_n (q;q^2)_{n} }
\\ =\sum_{n=0}^\infty q^{18n^2 + 3n}(1-q^{30n+15})  -q^3 \sum_{n=0}^\infty
q^{18n^2 + 15n}(1-q^{6n+3}) \label{ft8}
\end{multline}
\begin{equation}
\sum_{n=0}^\infty \frac{ (-1)^n q^{n(n+1)} (q^6;q^6)_n}{ (q^{2};q^2)_{2n+1} }
=\sum_{n=0}^\infty  q^{18n^2+6n} (1-q^{24n+12})   \label{ft9}
\end{equation}
\begin{multline}
\sum_{n=0}^\infty \frac{ (-1)^n q^{n(n+3)} (q^3;q^6)_n}{
(q^{2};q^4)_{n+1} (-q^2;q^2)_n (q;q^2)_{n} }
\\=\sum_{n=0}^\infty q^{18n^2+9n} (1-q^{18n+9})
 + q^2 \sum_{n=0}^\infty q^{18n^2+15n} (1-q^{6n+3}) \label{ft10}
\end{multline}

 In \S\ref{StdResults}, we will review some standard definitions and results to be used in the sequel.
In \S\ref{Proofs}, we indicate the Bailey pairs necessary to prove Identities~\eqref{m18-1}--\eqref{ft10}
and provide the keys to proving Identities~\eqref{m18-1}--\eqref{m24s-m5}.  In \S\ref{FT}, we will discuss and prove the false theta series identities~\eqref{ft1}--\eqref{ft10}.
Finally, in \S\ref{Lie} we discuss possible connections between
Identities~\eqref{m18-1}--\eqref{m18-4} and the standard level 6 modules associated with
the Lie algebra $A_{2}^{(2)}$.

\section{Standard definitions and results}\label{StdResults}

We will require a number of definitions and theorems from the
literature. It will be convenient to adopt Ramanujan's notation
for theta functions~\cite[p. 11, Eqs. (1.1.5)--(1.1.8)]{AB05}.
{\allowdisplaybreaks
\begin{defn} For $|ab|<1$, let
\begin{align} f(a,b) &:= \sum_{n=-\infty}^\infty a^{n(n+1)/2} b^{n(n-1)/2}, \label {fdef}\\
\varphi(q) &:=  f(q,q), \label{phidef}\\
\psi(q) &:= f(q,q^3), \label{psidef}\\
f(-q) &:= f(-q,-q^2).\label{PNSdef}
\end{align}
\end{defn}
}
 Both the Jacobi triple product identity and the quintuple
product identity were used extensively by Ramanujan (cf.~\cite{AB05}, \cite{AB07})
and Slater~\cite{S52}.  Rogers, on the other hand, appears to have
been unaware of the quintuple product identity, since he referred to~\cite[p. 333, Eq. (16)]{R94}
\begin{equation}
\label{remarkable}
\frac{(q^2;q^2)_\infty }{ (q^{30}; q^{30})_\infty (q;q^5)_\infty (q^4;q^5)_\infty} =
 (q^{13}; q^{30})_\infty (q^{17};q^{30})_\infty + q (q^7;q^{30})_\infty (q^{23};q^{30})_\infty,
\end{equation}
which follows immediately from the quintuple product identity,
as a ``remarkable identity" after observing that both sides of~\eqref{remarkable}
are equal to the same series.  Accordingly, we have chosen the name
``Ramanujan-Slater type identities" in our title for the identities in
this paper rather than ``Rogers-Ramanujan
type identities."

 Many proofs of the Jacobi triple product identity are known; see,
e.g.,~\cite[pp. 496--500]{AAR99} for two proofs. For a history and
many proofs of the quintuple product identity, see S. Cooper's
excellent survey article~\cite{C06}.
\begin{thm}[Jacobi's triple product identity] For $|ab|<1$,
   \begin{equation} \label{jtp}
     f(a,b) = (-a, -b, ab ; ab)_\infty.
   \end{equation}
\end{thm}

\begin{thm}[Quintuple product identity] For $|w|<1$ and $x\neq 0$,
  \begin{multline} \label{qpi}
    f(-wx^3, -w^2 x^{-3}) + x f(-wx^{-3}, -w^2 x^3) = \frac{ f(w/x, x) f(-w/x^2, -wx^2) }{ f(-w^2) } \\
      = (-wx^{-1}, -x, w; w)_\infty (wx^{-2}, wx^2; w^2)_\infty.
  \end{multline}
\end{thm}

The following is a special case of Bailey's ${}_6 \psi_6$ summation formula~\cite[Eq. (4.7)]{B36}
which appears in Slater~\cite[p. 464, Eq. (3.1)]{S51}.
\begin{thm}[Bailey]
\begin{multline} \label{6psi6}
 \sum_{r=-\infty}^\infty \frac{ (1-aq^{6r})(q^{-n};q)_{3r} (e;q^3)_r a^{2r} q^{3nr} }
 {(1-a)(aq^{n+1}; q)_{3r} (aq^3/e;q^3)_r e^r } \\
 = \frac{ (a;q^3)_\infty  (q^3/a;q^3)_\infty  (aq^2/e;q^3)_\infty  (aq/e; q^3)_\infty (q;q)_n (aq;q)_n (a^2/e; q^3)_n}
    { (q;q^3)_\infty (q^2;q^3)_\infty (q^3/e;q^3)_\infty (a^2/e; q^3)_\infty (a;q)_{2n} (aq/e;q)_n },
 \end{multline}
 where $a$ must be a power of $q$ so that the series terminates below.
\end{thm}

The next two $q$-hypergeometric summation formulas are due to
to Andrews~\cite[p. 526, Eqs. (1.8) and (1.9) respectively]{A73}.
\begin{thm}[$q$-analog of Gauss's ${}_2 F_{1} (\frac 12)$ sum]
\begin{equation} \label{q2ndGauss}
\sum_{n=0}^\infty \frac{ q^{n(n+1)} (a;q^2)_n (b;q^2)_n }
  { (q^2;q^2)_n (abq^2;q^4)_n } =
\frac{ (aq^2;q^4)_\infty (bq^2;q^4)_\infty }
{ (q^2;q^4)_\infty (abq^2;q^4)_\infty }.
\end{equation}
\end{thm}
{\allowdisplaybreaks
\begin{thm}[$q$-analog of Bailey's ${}_2 F_1 (\frac 12)$ sum]
\begin{equation} \label{qBailey}
 \sum_{n=0}^\infty \frac{ (bq;q^2)_n (b^{-1}q;q^2)_n c^n q^{n^2} }
{(cq;q^2)_n (q^4;q^4)_n }
= \frac{ (b^{-1}cq^2;q^4)_\infty (bcq^2; q^4)_\infty }{ (cq;q^2)_\infty }.
\end{equation}
\end{thm}
}
\begin{defn}
 A pair of sequences
\[
\left( \{\alpha_n (a,q) \}_{n=0}^\infty, \{
\beta_n(a,q)\}_{n=0}^\infty \right)\]
 is called
 a \emph{Bailey pair relative to $a$} if
 \begin{equation} \label{BPdef}
    \beta_n (a,q) = \sum_{r=0}^n \frac{\alpha_r (a,q)}{(a q;q )_{n+r} (q;q)_{n-r}}.
\end{equation}
\end{defn}
Bailey~\cite[p. 3, Eq. (3.1)]{B49} proved a key result, now known as ``Bailey's lemma,"
which led to the discovery of many Rogers-Ramanujan
type identities.

We will require several special cases of Bailey's lemma.
\begin{thm} If $\left( \{ \alpha_n (a,q) \}, \{ \beta_n (a,q) \} \right)$
form a Bailey pair, then {\allowdisplaybreaks
\begin{align}
\sum_{n=0}^\infty  a^n q^{n^2} \beta_n(a,q) &= \frac{1}{(aq;q)_\infty}
  \sum_{r=0}^\infty  a^r q^{r^2} \alpha_r (a,q) \label{aPBL} \\
\sum_{n=0}^\infty a^n q^{n^2} (-q;q^2)_n \beta_n(a,q^2) &= \frac{(-aq;q^2)_\infty }{(aq^2;q^2)_\infty}
  \sum_{r=0}^\infty a^r q^{r^2}  \alpha_r (a,q^2) \label{aTBL}\\
\frac{1}{1-q^2}\sum_{n=0}^\infty  q^{n(n+1)} (-q^2;q^2)_n \beta_n(q^2,q^2)
    &= \frac{1}{\s(-q^2)}
  \sum_{r=0}^\infty  q^{r(r+1)}  \alpha_r (q^2,q^2). \label{S2BL}  \\
\frac{1}{1-q^2}\sum_{n=0}^\infty  (-1)^n q^{n(n+1)} (q^2;q^2)_n \beta_n(q^2,q^2)
    &=
  \sum_{r=0}^\infty  (-1)^r q^{r(r+1)}  \alpha_r (q^2,q^2). \label{FBL}
 \end{align}
 }
\end{thm}

Eq.~\eqref{aPBL} is~\cite[p. 3, Eq. (3.1) with $\rho_1, \rho_2\to\infty$]{B49}.
Eq.~\eqref{aTBL} is~\cite[p. 3, Eq. (3.1) with $\rho_1=-\sqrt{q};\ \rho_2\to\infty$]{B49}.
Eq.~\eqref{S2BL} is~\cite[p. 3, Eq. (3.1) with $\rho_1=-q;\ \rho_2\to\infty$]{B49}.
Eq.~\eqref{FBL} is~\cite[p. 3, Eq. (3.1) with $\rho_1=q;\ \rho_2\to\infty$]{B49}.

\section{Proofs of Identities~\eqref{m18-1}--\eqref{m24s-m5}}\label{Proofs}

To facilitate the proofs of many of the identities, we will first need to establish a
number of Bailey pairs.  For instance,

\begin{lem} \label{BP2}
If {\allowdisplaybreaks
\begin{equation*}
  \alpha_n (1,q) =
   \left\{ \begin{array}{ll}
             1 &\mbox{if $n=0$}\\
             q^{\frac 92 r^2 - \frac 32 r} (1+q^{3r}) &\mbox{if $n=3r>0$}\\
             -q^{\frac 92 r^2 - \frac 92 r + 1} &\mbox{if $n=3r-1$}\\
             -q^{\frac 92 r^2 + \frac 92 r + 1} &\mbox{if $n=3r+1$}
              \end{array} \right.
\end{equation*} and
\[ \beta_n (1,q) = \frac  {(-1; q^3)_n } { (q;q)_{2n} (-1; q)_n },\]
then $\left( \alpha_n (1,q) , \beta_n (1,q) \right)$
form a Bailey pair relative to $1$.}
\end{lem}
\begin{pf}
Set $a=q$ and $e=-q^2$ in~\eqref{6psi6} and simplify to obtain
\begin{equation} \label{6psi6spec}
\sum_{r\in\mathbb Z} \frac{ (1-q^{6r+1}) q^{\frac 92 r^2 -\frac 32 r} }
 {(q;q)_{n-3r} (q; q)_{n+3r+1} }
 = \frac{   (-1; q^3)_n}
    {   (q;q)_{2n} (-1;q)_n }.
 \end{equation}
\begin{align*}
&\qquad\quad\sum_{r=0}^n \frac{ \alpha_r(1,q)}{(q;q)_{n-r} (q;q)_{n+r}} \\
&= \frac{1}{(q;q)_n^2} + \sum_{r\geqq 1} \frac{ \alpha_{3r}(1,q)}{(q;q)_{n-3r} (q;q)_{n+3r}}
         +  \sum_{r\geqq 1} \frac{ \alpha_{3r-1}(1,q)}{(q;q)_{n-3r+1} (q;q)_{n+3r-1}} \\
         &\qquad\qquad
         + \sum_{r\geqq 0} \frac{ \alpha_{3r+1}(1,q)}{(q;q)_{n-3r-1} (q;q)_{n-3r+1}} \\
 &=  \sum_{r\in\mathbb Z} \frac{ q^{\frac 92 r^2 -\frac 32 r} }{(q;q)_{n-3r} (q;q)_{n+3r}}
         -  \sum_{r\in\mathbb Z} \frac{ q^{\frac 92 r^2 + \frac 92 r + 1}}{(q;q)_{n+3r+1} (q;q)_{n-3r-1}}
         \\
 &=  \sum_{r\in\mathbb Z} \frac{ q^{\frac 92 r^2 -\frac 32 r} }{(q;q)_{n-3r} (q;q)_{n+3r+1}}
    \left( (1-q^{n+3r+1}) - q^{6r+1} (1-q^{n-3r}) \right) \\
&=\sum_{r\in\mathbb Z} \frac{ q^{\frac 92 r^2 -\frac 32 r} (1-q^{6r+1}) }{(q;q)_{n-3r} (q;q)_{n+3r+1}}
   = \frac{(-1;q^3)}{(q;q)_{2n} (-1;q)_n} \mbox{ (by~\eqref{6psi6spec}) }. \qed
\end{align*}
\end{pf}

The other necessary Bailey pairs can be established similarly, so we omit the details
and summarize the results in Table~\ref{BPtable}.1.

With the required Bailey pairs in hand, the identities can be proved.
For example, to prove Identity~\eqref{m18-1}, we proceed as follows:
\begin{pf}
Insert the Bailey pair P2 into Eq.~\eqref{aPBL} with $a=1$ to obtain
\begin{align*}
 & \quad\qquad \sum_{n=0}^\infty \frac{ q^{n(n+1)} (-1;q^3)_n}{ (q;q)_{2n} (-1;q)_n}\\ & =
\frac{1}{(q;q)_\infty}\left( 1 + \sum_{r=1}^\infty q^{\frac{27}{2} r^2 - \frac 32 r} (1+q^{3r})
 - \sum_{r=1}^\infty q^{\frac{27}{2} r^2 - \frac {15}{2} r +1}
  - \sum_{r=0}^\infty q^{\frac{27}{2} r^2 + \frac{15}{2} r+1 } \right) \\
& = \frac{1}{(q;q)_\infty} \left( \sum_{r=-\infty}^\infty q^{\frac{27}{2} r^2 - \frac 32 r}
  -q  \sum_{r=-\infty}^\infty q^{\frac{27}{2} r^2 - \frac {15}{2}r } \right)\\
  &= \frac{ f(q^{12}, q^{15}) - q  f(q^6, q^{21}) }{ f(-q)} \qquad \qquad\mbox{ (by~\eqref{jtp})} \\
  & =  \frac{ (q, q^8 , q^9; q^9)_\infty (q^7, q^{11};q^{18})_\infty }{ (q;q)_\infty}
  \qquad\mbox{  (by~\eqref{qpi}) }. \qed
\end{align*}
\end{pf}
The details of the proofs of the other identities are similar and therefore omitted, with the
key information summarized in Table~\ref{IdTable}.2.

\begin{landscape}
\begin{table}[hbt] \label{BPtable}
\centering

\begin{tabular}{|c|c|c|c|c|c|c|c|  }
\hline\hline
 & $a$ &$e$ & $\beta_n$ & $\alpha_{3r+1}$ & $\alpha_{3r}$ & $\alpha_{3r-1}$ & rel to \\
\hline
P1& $q$ & $-q^2$ & $\frac{(-1;q^3)}{(q;q)_{2n} (-1;q)_n}$
      &  $-q^{\frac 92 r^2 + \frac 92 r + 1} $
      & $ q^{\frac 92 r^2 - \frac 32 r} (1+q^{3r})$
      & $ -q^{\frac 92 r^2 - \frac 92 r + 1} $
      & $1$ \\
\hline
P2& $q$ & $-q^2$ &$\frac  {q^n (-1; q^3)_n } { (q;q)_{2n} (-1; q)_n }$
        &  $-q^{\frac 92 r^2 + \frac 32 r}$
        &  $q^{\frac 92 r^2 - \frac 32 r} (1+q^{3r}) $
        &  $-q^{\frac 92 r^2 - \frac 32 r}$
        & $1$ \\
\hline
P3& $q^2$ & $-q$ & $ \frac  {(-q^3; q^3)_n } { (q^2;q)_{2n} (-q; q)_n }$
               & $-2q^{\frac 92 r^2 + \frac 92 r + 1} $
               & $ q^{\frac 92 r^2 + \frac 32 r}$
                & $q^{\frac 92 r^2 - \frac 32 r} $
                & $q$\\
 \hline
P4&    $q^2$  & $-q^{5/2} $ & $\frac  {(-q^{3/2}; q^3)_n } { (q^2;q)_{2n} (-q^{1/2}; q)_n }$
         & $-q^{\frac 92 r^2+3r+\frac 12} (1+q^{3r+\frac 32})$
          & $q^{\frac 92 r^2}$
          & $q^{\frac 92 r^2}$
          & $q$ \\
 \hline
P5& $q^2$ & $-q^{5/2}$ & $ \frac  {q^n (-q^{3/2}; q^3)_n } { (q^2;q)_{2n} (-q^{1/2}; q)_n }$
     & $ -q^{\frac 92 r^2} (q^{6r+\frac 32}+q^{3r})$
     & $q^{\frac 92 r^2+3r} $
     & $q^{\frac 92 r^2-3r} $
     & $q$\\
 \hline
P6 & $q$ & $-q^2$ & $ \frac  {(1-q)(-1; q^3)_n } { (q;q)_{2n} (-1; q)_n }$
      & $0$
      & $ q^{\frac 92 r^2 - \frac 32 r} (1-q^{6r+1}) $
      & $ -q^{\frac 92 r^2 - \frac 92 r + 1}(1-q^{6r-1})$
      & $q$\\
 \hline
 P7 & $q$ & $q^2$
       & $\frac  {(q^3; q^3)_{n-1} } { (q^2;q)_{2n-1} (q;q)_{n-1} }$
       & $(-1)^{r+1} q^{\frac 92 r^2 + \frac 32 r +1}\frac{1-q^{6r+3}}{1-q}$
       & $(-1)^r  q^{\frac 92 r^2 - \frac 32 r} \frac{1-q^{6r+1}}{1-q} $
       & $(-1)^{r+1}  q^{\frac 92 r^2 - \frac 92 r+1} \frac{1-q^{6r-1}}{1-q}$
       & $q$\\
 \hline
\end{tabular}
\caption{By specializing $a$ and $e$ in~\eqref{6psi6} as indicated, each of the following
seven Bailey pairs (relative to $1$ or $q$ as stated) can be established.
In all cases $\alpha_0 = \beta_0 = 1$.}
\end{table}
\end{landscape}

\begin{table} \label{IdTable}
\caption{Proofs of identities~\eqref{m18-1}--~\eqref{m24s-m5}}
\begin{tabular}{|c|c|c|c|c| }
\hline\hline
Eq. & Bailey& Bailey  & $a$ &  \\
 &  pair &  lemma &  &  \\
\hline
\eqref{m18-1} & P2 & \eqref{aPBL} & $1$ &\\
\eqref{m18-2} & P1 & \eqref{aPBL} & $1$&\\
\eqref{m18-3} & P3 & \eqref{aPBL} & $q$&\\
\eqref{m18-4} &  $-$  & $-$  & $-$ & $q^{-1} \times ( \eqref{m18-2} - \eqref{m18-1} ) $\\
\eqref{m18-m1} & $-$ & $-$ & $-$ & ~\cite[p. 433, (B4) $+q\times$(B2)]{B47}\\
\eqref{m18-m2} &  $-$ & $-$ & $-$ & ~\cite[p. 433, (B4)$+q^2\times$(B1)]{B47}\\
\eqref{m18-m3} & $-$ & $-$ & $-$  &~\cite[p. 433, (B3)]{B47}\\
\eqref{m18-m4} & $-$ & $-$ & $-$  &~\cite[p. 433, (B2) $-q\times$(B1)]{B47}\\
\eqref{m24t-2} & $-$ & $-$ & $-$  &  Set $b=e^{\pi i/3}$ and $c=1$ in \eqref{qBailey}.\\
\eqref{m24t-1} & P2 & \eqref{aTBL} & $1$ &\\
\eqref{m24t-3} & P1 & \eqref{aTBL} & $1$ &\\
\eqref{m24t-4} & $-$ & $-$ & $-$  &Set $b=e^{\pi i/3}$ and $c=q^2$ in~\eqref{qBailey}.\\
\eqref{m24t-5} & $-$ & $-$ & $-$  & $q^{-1}\times (\eqref{m24t-3}-\eqref{m24t-1})$\\
\eqref{m24t-m2} & $-$ & $-$ & $-$  &Set $b=e^{2\pi i/3}$ and $c=1$ in \eqref{qBailey}.\\
\eqref{m24t-m1} & $-$ & $-$ & $-$  & \cite[p. 434, (C3)$ +q\times$(C2)]{B47}\\
\eqref{m24t-m3} & $-$ & $-$ & $-$  & \cite[p. 434, (C3)$ +q^3\times$(C1)]{B47}\\
\eqref{m24t-m4} & $-$ & $-$ & $-$  & Set $b=e^{2\pi i/3}$ and $c=q^2$ in \eqref{qBailey}\\
\eqref{m24t-m5} & $-$ & $-$ & $-$  & $q^{-1}\times (\eqref{m24t-m1}-\eqref{m24t-m3})$\\
\eqref{m24s-1} & $-$ & $-$ & $-$  & $\eqref{m24s-3}-q\times\eqref{m24s-5} $ \\
\eqref{m24s-2} & $-$ & $-$ & $-$  &Set $a= e^{\pi i /3} $, $b=e^{-\pi i/ 3}$ in~\eqref{q2ndGauss}.\\
\eqref{m24s-3} & P4 & \eqref{S2BL} & $q$ &\\
\eqref{m24s-4} & $-$ & $-$ & $-$  &Set $a= e^{\pi i /3} q^2$, $b=e^{-\pi i/ 3} q^2$ in~\eqref{q2ndGauss}.\\
\eqref{m24s-5} & P5 & \eqref{S2BL} & $q$ &  \\
\eqref{m24s-m1} & $-$ & $-$ & $-$  & $\eqref{m24s-m3}+q\times\eqref{m24s-m5} $ \\
\eqref{m24s-m2} & $-$ & $-$ & $-$  & Set $a= e^{2\pi i /3} $, $b=e^{-2\pi i/ 3}$ in~\eqref{q2ndGauss}.\\
\eqref{m24s-m3} & J4~\cite[p. 149]{S52} & \eqref{S2BL} & $q$ &\\
\eqref{m24s-m4} &$-$ & $-$ & $-$  & \small{Set $a= e^{2\pi i /3} q^2$, $b=e^{-2\pi i/ 3} q^2$ in~\eqref{q2ndGauss}.}\\
\eqref{m24s-m5} &J5~\cite[p. 149]{S52} & \eqref{S2BL} & $q$ & \\
\hline
\end{tabular}
\end{table}

\section{False theta series identities}\label{FT}
Rogers introduced the term ``false theta series" and included a number of
related identities in his 1917 paper~\cite{R17}.  Ramanujan presented a number
of identities involving false theta series in his lost notebook~\cite[p. 256--259, \S11.5]{AB05}.

 Recalling that Ramanujan defines the theta function as
  \begin{align*} f(a,b)&:= \sum_{n=-\infty}^\infty a^{n(n+1)/2} b^{n(n-1)/2}\\
  &=
   \sum_{n=0}^\infty a^{n(n+1)/2} b^{n(n-1)/2} + \sum_{n=1}^\infty a^{n(n-1)/2} b^{n(n+1)/2}\\
 &= 1 +a + b + a^3 b + ab^3 + a^6 b^3 + a^3 b^6 + a^{10} b^6 + a^6 b^{10} + \dots,
  \end{align*}
 let us define the corresponding \emph{false theta function} as
\begin{align*}
\ft(a,b)&:=\sum_{n=0}^\infty a^{n(n+1)/2} b^{n(n-1)/2} - \sum_{n=1}^\infty a^{n(n-1)/2} b^{n(n+1)/2}\\
&= \sum_{n=0}^\infty a^{n(n+1)/2} b^{n(n-1)/2} (1 - b^{2n+1}) \\
&=1 +a - b + a^3 b - ab^3 + a^6 b^3 - a^3 b^6 + a^{10} b^6 - a^6 b^{10} + \dots.
\end{align*}
In practice, $a$ and $b$ are always taken to be $\pm q^h$ for some integer or half-integer $h$.

The key to the proof of each false theta series identity is indicated in Table~\ref{FTtable}.1.

\begin{table} \label{FTtable}
\caption{Proofs of identities~\eqref{ft1}--~\eqref{ft10}}
\begin{tabular}{|c|c|c|c|c| }
\hline\hline
Eq. & Bailey pair & form of Bailey lemma & $a$ &  \\
\hline
\eqref{ft1} & $-$ & $-$ & $-$ & \eqref{ft3}$-q\times$\eqref{ft5}\\
\eqref{ft2} & P6 & \eqref{FBL} &$q$ & \\
\eqref{ft3} & P4 & \eqref{FBL} &$q$ &\\
\eqref{ft4} & P3 & \eqref{FBL} &$q$ &\\
\eqref{ft5} & P5 & \eqref{FBL} &$q$ &\\
\eqref{ft6} & $-$ & $-$ & $-$ & \eqref{ft8}+$q\times$\eqref{ft10}\\
\eqref{ft7} & P7 & \eqref{FBL} &$q$ &\\
\eqref{ft8} & J4~\cite[p. 149]{S52} & \eqref{FBL} &$q$ &\\
\eqref{ft9} & $-$ & $-$ & $-$ & See~\cite[Entry 5.4.2]{AB07}\\
\eqref{ft10} & J5~\cite[p. 149]{S52} & \eqref{FBL} &$q$ &\\
\hline
\end{tabular}
\end{table}

\section{Connections with Lie algebras}\label{Lie}

Let $\g$ be the affine Kac-Moody Lie algebra $A_1^{(1)}$ or $A_2^{(2)}$.
Let $h_0, h_1$ be the usual basis of a maximal toral
subalgebra $T$ of $\g$.
Let $d$ denote the ``degree derivation" of $\g$
and $\tilde{T}:= T \oplus \mathbb C d$.
For all dominant integral $\lambda\in\tilde{T}^*$, there is an
essentially
unique irreducible, integrable, highest weight module $L(\lambda)$,
assuming without loss of generality that $\lambda(d) = 0$.
Now $\lambda= s_0 \Lambda_0 + s_1 \Lambda_1$ where
$\Lambda_0$ and $\Lambda_1$ are the fundamental weights,
given by $\Lambda_i(h_j) = \delta_{ij}$ and $\Lambda_i(d) = 0$;
here $s_0$ and $s_1$ are nonnegative integers.
For $A_1^{(1)}$, the canonical central element is
$c= h_0 + h_1$, while for
$A_2^{(2)}$, the canonical central element is
$c = h_0 + 2h_1$.
The quantity $\lambda(c)$ (which equals $s_0+s_1$ for $A_1^{(1)}$ and
which equals $s_0+2s_1$ for $A_2^{(2)}$)
is called the \emph{level} of $L(\lambda)$. (cf.\cite{K90}, \cite{LM78}.)

Additionally (see~\cite{LM78}),  there is an
infinite product $F_{\g}$ associated with $\g$, often light-heartedly
called the ``fudge factor," which needs to be divided out of the
the principally specialized character
$\chi(L(\lambda))
= \chi(s_0 \Lambda_0 + s_1\Lambda_1)$, in order to obtain
the quantities of interest here.  For $\g=A_1^{(1)}$,
the fudge
factor is given by
$F_{\g} = (q;q^2)_\infty^{-1}$, while for $\g=A_2^{(2)}$, it is
given by
$F_{\g} = \left[ (q;q^6)_\infty (q^5;q^6)_\infty \right]^{-1}$.

 Now $\g$ has a certain infinite-dimensional Heisenberg subalgebra
known as the ``principal Heisenberg vacuum subalgebra" $\sfrak$
(see~\cite{LW78} for the construction of $A_1^{(1)}$ and~\cite{KKLW81}
for that of $A_2^{(2)}$).
As shown in~\cite{LW82}, the principal character
$\chi(\Omega(s_0 \Lambda_0 + s_1 \Lambda_1))$, where $\Omega(\lambda)$
is the vacuum space for $\sfrak$ in $L(\lambda)$, is
\begin{equation} \label{char}
\chi(\Omega(s_0 \Lambda_0 + s_1\Lambda_1))
=  \frac{\chi( L(s_0 \Lambda_0 + s_1 \Lambda_1)) }{ F_{\g}  },
\end{equation}
where $\chi(L(\lambda))$ is the principally specialized character of
$L(\lambda)$.

 By~\cite{LM78} applied to~\eqref{char} in the case of $A_1^{(1)}$,
for standard modules of odd level $2k+1$,
$$\chi(\Omega( (2k-i+2)\Lambda_0 + (i-1)\Lambda_1 ))$$
is given by
Andrews' analytic generalization of the Rogers-Ramanujan
identities~\cite{A74}:
\begin{equation}\label{AndGor}
\sum_{n_1, n_2, \dots, n_{k}\geqq 0}
 \frac{ q^{N_1^2 + N_2^2 + \cdots + N_{k}^2 + N_i+N_{i+1}+\cdots+N_{k}}}
 {(q;q)_{n_1} (q;q)_{n_2} \cdots (q;q)_{n_{k}}  }
= \frac{(q^i,q^{2k+3-i},q^{2k+3};q^{2k+3})_\infty }{(q;q)_\infty },
\end{equation}
where $1\leqq i \leqq k+1$ and $N_j: = n_j + n_{j+1} + \cdots + n_{k}$.
The combinatorial counterpart to~\eqref{AndGor} is Gordon's
partition theoretic generalization of the Rogers-Ramanujan
identities~\cite{G61}; this generalization was explained
vertex-operator theoretically in~\cite{LW84} and~\cite{LW85}.

In addition, for the $A_1^{(1)}$ standard modules of even level $2k$,
\[ \chi( \Omega( (2k-i+1)\Lambda_0 +  (i-1)\Lambda_1 )) \]
is given
by Bressoud's
analytic identity~\cite[p. 15, Eq. (3.4)]{B80}
\begin{equation}\label{BressoudEven}
\sum_{n_1, n_2, \dots, n_{k}\geqq 0}
 \frac{ q^{N_1^2 + N_2^2 + \cdots + N_{k}^2 + N_i+N_{i+1}+\cdots+N_{k}}}
 {(q;q)_{n_1} (q;q)_{n_2} \cdots (q;q)_{n_{k-1}} (q^2;q^2)_{n_k}  }
= \frac{(q^i,q^{2k+2-i},q^{2k+2};q^{2k+2})_\infty }{(q;q)_\infty },
\end{equation}
where $1\leqq i \leqq k+1$,
and its partition theoretic counterpart~\cite[p. 64, Theorem, $j=0$ case]{B79};
likewise, this generalization was explained
vertex-operator theoretically in~\cite{LW84} and~\cite{LW85}.

  Notice that the infinite products associated with
level $\ell$ standard modules for $A_1^{(1)}$ in~\eqref{AndGor}
and~\eqref{BressoudEven}
are instances of the Jacobi triple product
identity for modulus $\ell+2$ divided by $(q;q)_\infty$.

   Probably the most efficient way of deriving~\eqref{AndGor} is via
the Bailey lattice~\cite{AAB87}, which is an extension of the Bailey chain
concept (\cite{A84}; cf. \cite[\S 3.5, pp. 27ff]{A86})
built upon the ``unit Bailey pair"
\[ \beta_n(1,q) = \left\{
   \begin{array}{ll}
      1 &\mbox{if $n=0$}\\
      0 &\mbox{if $n>0$}
     \end{array}  \right. \]
\[ \alpha_n(1,q) = \left\{
     \begin{array}{ll}
        1 &\mbox{if $n=0$}\\
        (-1)^n q^{n(n-1)/2} (1+q^n) &\mbox{if $n>0$.}
     \end{array} \right. \]

  Similarly,
\eqref{BressoudEven} follows from a Bailey lattice built upon the
Bailey pair

\[ \beta_n(1,q) = \frac{1}{(q^2;q^2)_n}, \]
\[ \alpha_n(1,q) = \left\{
     \begin{array}{ll}
        1 &\mbox{if $n=0$}\\
        (-1)^n 2 q^{n^2} &\mbox{if $n>0$.}
     \end{array} \right. \]

Thus the standard modules of $A_1^{(1)}$ may be compactly ``explained"
via two interlaced instances of the Bailey lattice.

  In contrast, the standard modules of $A_{2}^{(2)}$ are not as well
understood, and a uniform $q$-series and partition correspondence
analogous to what is known for $A_1^{(1)}$ has thus far remained
elusive.

As with $A_1^{(1)}$, there are $1+\lfloor \frac{\ell}{2}
\rfloor$ inequivalent level $\ell$ standard modules associated with
the Lie algebra $A_2^{(2)}$, but the analogous
quantity for the level $\ell$ standard modules
\[ \chi (\Omega( (\ell-2i+2)\Lambda_0 + (i-1)\Lambda_1 )) \]
is given
by instances of the quintuple product identity
(rather than the triple product identity) divided by $(q;q)_\infty$:

\begin{equation} \label{A22prodside}
 \frac{ (q^i, q^{\ell+3-i}, q^{\ell+3}; q^{\ell+3})_\infty (q^{\ell+3-2i},
  q^{\ell+2i+3}; q^{2\ell+6})_\infty}{(q;q)_\infty},
\end{equation}
where $1\leqq i \leqq 1 + \lfloor \frac{\ell}{2} \rfloor $; see~\cite{LM78}.

  It seems quite plausible that in the case of $A_2^{(2)}$, the analog of
the Andrews-Gordon-Bressoud identities would involve the interlacing of
six Bailey lattices in contrast to the two that were necessary for
$A_1^{(1)}$.   To see this, consider the following set of Andrews-Gordon-Bressoud type
identities where the product sides involve instances of
the quintuple product identity
rather than the triple product identity:

{\allowdisplaybreaks
\begin{multline} \label{lev6k2}
 \sum_{n_1, n_2, \dots, n_k \geqq 0}
  \frac{ q^{ N_1(N_1+1)/2 + N_2(N_2+1) + N_3(N_3+1)+ \cdots +N_k(N_k+1) + N_k^2}  }
  { (q;q)_{n_1} (q;q)_{n_2} \cdots (q;q)_{n_{k-1}} (q;q)_{2n_k+1} (-q^{N_1+1};q)_\infty}
\\= \frac{ (q^{k}, q^{5k-1}, q^{6k-1}; q^{6k-1})_\infty (q^{4k-1},q^{8k-1}; q^{12k-2})_\infty }{(q;q)_\infty}
\end{multline}
\begin{multline} \label{lev6k3}
\sum_{n_1, n_2, \dots, n_{k+1} \geqq 0}
  \frac{ q^{N_1^2 + N_2^2 + \cdots +N_{k}^2} \left( \frac{n_k- n_{k+1}+1}{3} \right) }
  { (q;q)_{n_1} (q;q)_{n_2} \cdots (q;q)_{n_{k+1}} (q;q)_{2n_k- n_{k+1} }}
\\= \frac{ (q^{k}, q^{5k}, q^{6k}; q^{6k})_\infty (q^{4k},q^{8k}; q^{12k})_\infty }{(q;q)_\infty}
\end{multline}
\begin{multline} \label{lev6k4}
 \sum_{n_1, n_2, \dots, n_k \geqq 0}
  \frac{ q^{ N_1(N_1+1)/2 + N_2(N_2+1) + N_3(N_3+1)+ \cdots +N_k(N_k+1)}  }
  { (q;q)_{n_1} (q;q)_{n_2} \cdots (q;q)_{n_{k-1}} (q;q)_{2n_k+1} (-q^{N_1+1};q)_\infty}
\\= \frac{ (q^{2k}, q^{4k+1}, q^{6k+1}; q^{6k+1})_\infty (q^{2k+1},q^{10k+1}; q^{12k+2})_\infty }{(q;q)_\infty}
\end{multline}
} {\allowdisplaybreaks
\begin{multline} \label{lev6k5}
\sum_{n_1, n_2, \dots, n_k \geqq 0}
  \frac{ q^{N_1^2 + N_2^2 + \cdots +N_{k-1}^2+2N_k^2}}
  { (q;q)_{n_1} (q;q)_{n_2} \cdots (q;q)_{n_{k-1}} (q;q)_{2n_k}}
\\= \frac{ (q^{k}, q^{5k+2}, q^{6k+2}; q^{6k+2})_\infty (q^{4k+2},q^{8k+2}; q^{12k+4})_\infty }{(q;q)_\infty}
\end{multline}
} {\allowdisplaybreaks
\begin{multline} \label{lev6k6}
\sum_{n_1, n_2, \dots, n_k \geqq 0}
  \frac{ q^{N_1^2 + N_2^2 + \cdots +N_k^2} (-1;q^3)_{n_k}}{ (q;q)_{n_1} (q;q)_{n_2} \cdots (q;q)_{n_{k-1}} (q;q)_{2n_k} (-1;q)_{n_k} }
\\= \frac{ (q^{k+1}, q^{5k+2}, q^{6k+3}; q^{6k+3})_\infty (q^{4k+1},q^{8k+5}; q^{12k+6})_\infty }{(q;q)_\infty}
\end{multline}
} {\allowdisplaybreaks
\begin{multline}\label{lev6k7}
\sum_{n_1, n_2, \dots, n_k \geqq 0}
  \frac{ q^{N_1^2 + N_2^2 + \cdots +N_k^2}}{ (q;q)_{n_1} (q;q)_{n_2} \cdots (q;q)_{n_{k-1}} (q;q)_{2n_k}}
\\= \frac{ (q^{k+1}, q^{5k+3}, q^{6k+4}; q^{6k+4})_\infty (q^{4k+2},q^{8k+6}; q^{12k+8})_\infty }{(q;q)_\infty},
\end{multline}
}
where
$\left( \frac{n}{p} \right)$ in~\eqref{lev6k3} is the Legendre symbol.
We note that~\eqref{lev6k3} first appeared in~\cite[p. 400, Eq. (1.7)]{S04} and that
\eqref{lev6k7} is due to Andrews~\cite[p. 269, Eq. (1.8)]{A84}.  While~\eqref{lev6k2},
\eqref{lev6k4}, and \eqref{lev6k5}
probably have not appeared explicitly in the literature, they each follow from
building a Bailey chain on a known Bailey pair and may be regarded as nothing
more than a standard exercise in light of Andrews' discovery of
the Bailey chain~(\cite{A84}; cf.~\cite[\S3.5]{A86}).
Indeed the $k=1$ cases of~\eqref{lev6k2},~\eqref{lev6k4},~\eqref{lev6k5}, and~\eqref{lev6k7}
are all due to Rogers and appear in Slater's list~\cite{S52} as Eqs. (62), (80), (83),
and (98) respectively.
On the other hand,~\eqref{lev6k6} is new since it arises from inserting
a new Bailey pair, namely the one from Lemma~\ref{BP2} in this paper, into the
Bailey chain mechanism.
Notice that as $k$ runs through the positive integers
in the numerators of the right hand sides of~\eqref{lev6k2}--\eqref{lev6k7},
we obtain instances of the quintuple product identity for all moduli represented
in~\eqref{A22prodside}
(except for the trivial level 1 case where the relevant identity
reduces to ``$1=1$").   It is because of the preceding observations
that we conjecture that $A_2^{(2)}$ may
be ``explained" by six interlaced Bailey lattices.

We now turn our attention to combinatorial considerations in the
context of $A_2^{(2)}$. In his 1988 Ph.D. thesis S.
Capparelli~\cite{C88} conjectured two beautiful partition identities
resulting from his analysis of the two inequivalent level 3 standard
modules of $A_2^{(2)}$, using the theory in~\cite{LW84}
and~\cite{LW85}.  Capparelli's conjectures were first proved by
Andrews~\cite{A94} using combinatorial methods.  Later, Lie
algebraic proofs were found by Tamba and Xie~\cite{TX95} and
Capparelli himself~\cite{C96}. More recently, Capparelli~\cite{C04}
related the principal characters of the vacuum spaces for the
standard modules of $A_2^{(2)}$ for levels 5 and 7 to some known
$q$-series and partition identities. In the same way, our
identities~\eqref{m18-1}--\eqref{m18-4} appear to correspond to the
standard modules for level 6.

\section*{Acknowledgements}
Many thanks are due to Jim Lepowsky and Robert Wilson for their help
with the exposition in \S\ref{Lie}. We also thank George Andrews for
his encouragement and several useful suggestions.  Finally, we thank
the referee for helpful comments.

\end{document}